\documentclass[12pt]{article}

\usepackage{amssymb}

      \def\hbfT{{ \hat T}}
      \def\hbfA4{{ A_{ijkl} (\hat T)}}
      \def\traT { tr ~\hat T}
      \def\trT2{ tr ~{\hat T}^2}

\begin{document}

\vskip0.2in
\centerline {\bf Isotropic tensor-valued polynomial functions of fourth-order tensors}

\vskip0.8in

\centerline {B. A.~Younis \footnote{Department of Civil \& Environmental Engineering,
University of California, Davis, CA 95616, U.S.A.}\& G. F. Smith \footnote{Department of Mechanical Engineering \& Mechanics,
 Lehigh University, Bethlehem, PA 18017--5120, USA}}

\vskip0.6in

\centerline{\bf Abstract}
\vskip0.3in
\noindent Fourth-order tensor-valued functions appear in numerous fields of study.  
The formulation of practical models for these complex functions often requires their representation 
 in terms of tensors of order two. In this paper, we develop
an appropriate representation formula by assuming that  the isotropic fourth--order tensor--valued
function is a polynomial function in the components of two symmetric second--order
tensors of degree $\leq$ 2. We illustrate the utility of the result by applying it to obtain a representation
of the fluctuating velocity, pressure-gradient correlations of turbulence.


\section{Introduction}

Fourth-order tensor-valued functions are frequently encountered in applied sciences e.g. in computational mechanics 
where they play an important role in the formulation of constitutive relations \cite{ITS}. In the study of turbulence, 
the importance of these functions becomes apparent when considering the
state of the art in the closure of the Reynolds--averaged form of the
Navier--Stokes equations. This level of closure requires the
solution of a differential transport equation for each  non--zero
component of the Reynolds--stress tensor $\overline {u_iu_j}$ (see, e.g., \cite{LRR},\cite{SSG}). The equations
governing the evolution of this tensor can be derived directly from the 
Navier--Stokes equations. For an incompressible fluid of uniform
properties,  they may be written as:
\begin{eqnarray}
\frac{\partial \overline{u_i u_j}}{\partial t} 
+ U_k \frac{\partial  \overline{u_{i} u_{j} }} {\partial x_k} 
= & - & \left ( \overline{u_{i} u_{k} }
\frac{ \partial U_{j} } {\partial x_{k}}
+ \overline{  u_{j} u_{k} } \frac{  \partial U_{i} }
{ \partial x_{k}} \right )  
  -    2 \nu \left(\overline{ \frac{  \partial u_{i}} { \partial x_{k}}
\frac{  \partial u_{j}}{ \partial x_{k} }} \right) 
\nonumber \\[3 mm]
 & - &  \frac{ \partial }{  \partial x_{k}}
 \left(\overline{  u_{i} u_{j} u_{k}}
 - \nu \frac{  \overline{  \partial u_{i} u_{j}}}
 {\partial x_{k}}\right) \nonumber \\[4 mm]
 & - & \frac{1}{\rho} \left( \overline{u_j
 \frac{  \partial p}{ \partial x_i}} +  \overline {u_i \frac{  \partial p}
{ \partial x_j}} \right)
\label{e:rs}
\end{eqnarray}       

\noindent where $U_i$ is the time--averaged velocity, $\rho$ and $\nu$ are the fluid
density and kinematic viscosity and $u_i$ and $p$ are the 
fluctuating velocity and pressure respectively. Repeated subscripts indicate
summation.

An exact expression the last term in Eq. (\ref{e:rs}) was  obtained by ~\cite{Chou}
by taking the divergence
of the Navier--Stokes equations to obtain a Poisson equation for the 
fluctuating pressure. The  general solution of this equation,
after multiplication  by 
the fluctuating velocity, is given as:                                          
\begin{eqnarray}
\frac{1}{\rho}
\overline{u_i \frac{  \partial p}{ \partial x_j}} & = & \frac{1}{2 \pi}
\int \int \int \left ( \frac{  \partial U'_k}{ \partial x_l}
\frac{   \partial \overline{ u'_{l} u_{i}}} {\partial x_{k}}
\right )_{,j}  \: \frac{1}{{\bf r}} \: d Vol   \nonumber \\
               & + & \frac{1}{4 \pi}
\int \int \int \left (\overline{ u'_{k} u'_{l} u_{i}} \right )_{,klj} \: 
\frac{1}{{\bf r}} \: d Vol
\label{e:sol}
\end{eqnarray}      
where the integrations extend over the whole moving fluid,  d Vol is
the volume element and terms with 
a prime relate to values at point which ranges over the region of the moving
fluid separated by distance ${\bf r}$ from the point  where 
the pressure fluctuations are evaluated. 

With the assumption of local homogeneity which is commonly invoked in turbulence studies, and by retaining only the part
 which gives rise to the fourth-order tensor, 
there results:
\begin{eqnarray}
\frac{1}{\rho} \left ( \overline{u_j
 \frac{  \partial p}{ \partial x_i}} +  \overline {u_i \frac{  \partial p}
{ \partial x_j}} \right )= (A_{ijkl}+A_{jikl})
\frac{  \partial U_k}{ \partial x_l}
\label{e:fin}
\end{eqnarray}
\noindent where the tensor functions $A_{ijkl}$,  
$A_{jikl}$ (by permutation of i and j) are defined by
Eq. (\ref{e:sol}).\\

The ability to use Eq. (\ref{e:rs}) to model the effects of turbulence on the motion of a fluid thus depends on the
ability to represent the 
fourth--order tensor functions $A_{ijkl}$ and $A_{jikl}$ in terms of second-order tensors. The purpose of this paper is to derive 
an appropriate representation for these functions, and to demonstrate the use of the result in this particular application.

\section{Isotropic fourth-order tensor-valued functions}
For convenience, we write:
\begin{eqnarray}
A_{ijkl} & = & \hbfA4, \qquad
\hat {T} = T_{ij}, \qquad T_{ij} = \overline{u_i u_j}
\label{e:1.2}
\end{eqnarray}

We proceed by assuming that  $\hbfA4$ is an isotropic fourth--order tensor--valued
polynomial function of the components $T_{ij}$ of the symmetric second--order
tensor $\hbfT$ of degree $\leq$ 2. Thus, we have
\begin{equation}
A_{ijkl} (\hat T) = \alpha_{ijklmn} T_{mn} + \beta_{ijklmnpq} T_{mn} T_{pq}.
\end{equation}

The tensors $\alpha_{ijklmn}$ and $\beta_{ijklmnpq}$ are isotropic tensors and are
expressible in terms of outer products of the Kronecker delta tensor $\delta_{ij}$.

The sixth--order isotropic tensor $\alpha_{ijklmn}$ is expressible as a linear
combination of 15 distinct isomers of $\delta_{ij}$ $\delta_{kl}$ $\delta_{mn}$
\cite{Smith70}. 
Tensors arising from $\delta_{ij}$ $\delta_{kl}$ $\delta_{mn}$ upon permuting the 
subscripts i j $\ldots$ n are  referred to as isomers of 
$\delta_{ij}$ $\delta_{kl}$ $\delta_{mn}$. We have:
\begin{eqnarray}
\alpha_{ijklmn} & = & \alpha_1 \delta_{ij} \delta_{kl} \delta_{mn}  
                    + \alpha_2 \delta_{ij} \delta_{km} \delta_{ln} 
                    + \alpha_3 \delta_{ij} \delta_{kn} \delta_{lm} 
                       \nonumber \\[3 mm]
                & + & \alpha_4  \delta_{ik} \delta_{jl} \delta_{mn}
                    + \alpha_5  \delta_{ik} \delta_{jm} \delta_{ln}
                    + \alpha_6  \delta_{ik} \delta_{jn} \delta_{lm}
                       \nonumber \\[3 mm]
                & + & \alpha_7  \delta_{il} \delta_{jk} \delta_{mn}
                    + \alpha_8  \delta_{il} \delta_{jm} \delta_{kn}
                    + \alpha_9  \delta_{il} \delta_{jn} \delta_{km}
                        \nonumber \\[3 mm] 
                & + & \alpha_{10} \delta_{im} \delta_{jk} \delta_{ln}
                    + \alpha_{11} \delta_{im} \delta_{jl} \delta_{kn}
                    + \alpha_{12} \delta_{im} \delta_{jn} \delta_{kl}
                         \nonumber \\[3 mm] 
                & + & \alpha_{13} \delta_{in} \delta_{jk} \delta_{lm}
                    + \alpha_{14} \delta_{in} \delta_{jl} \delta_{km}
                    + \alpha_{15} \delta_{in} \delta_{jm} \delta_{kl}
\end{eqnarray}
The eighth--order isotropic tensor $\beta_{ijklmnpq}$ is expressible as a linear
combination of 105 distinct isomers of $\delta_{ij}$ $\delta_{kl}$ $\delta_{mn}$
$\delta_{pq}$. We thus have
\begin{eqnarray}
A_{ijkl} & =  &(\alpha_{1} \delta_{ij} \delta_{kl} \delta_{mn} + \ldots + 
            \alpha_{15} \delta_{in} \delta_{jm} \delta_{kl}) T_{mn} \nonumber \\
         & + & (\beta_{1} \delta_{ij} \delta_{kl} \delta_{mn} \delta_{pq} + \ldots +
             \beta_{105} \delta_{iq} \delta_{jp} \delta_{kn} \delta_{lm}) T_{mn} T_{pq} 
\label{e:1.8}
\end{eqnarray}

Only 91 of the 105 distinct isomers of 
$\delta_{ij}$ $\delta_{kl}$ $\delta_{mn}$ $\delta_{pq}$
are linearly independent. Consider tensors of the form:
\begin{eqnarray}
\delta^{i k m p}_{j l n q} =   \left | \begin{array} {cccc}
\delta_{ij} & \delta_{il} & \delta_{in} & \delta_{iq} \\ 
\delta_{kj} & \delta_{kl} & \delta_{kn} & \delta_{kq} \\ 
\delta_{mj} & \delta_{ml} & \delta_{mn} & \delta_{mq} \\ 
\delta_{pj} & \delta_{pl} & \delta_{pn} & \delta_{pq} 
\end{array}
 \right |.
\label{e:1.9}
\end{eqnarray}

If the tensor is three--dimensional, it is a null tensor. For any of the $3^8$ possible
choices of values which i, $\dots$, q may assume, at least two rows of the
determinant will be the same and the component will be zero. With (\ref{e:1.9}), we have:
\begin{eqnarray}
\delta^{i k m p}_{j l n q} T_{mn} T_{pq} & = & 2 (T_{ij} T_{kl} - T_{il} T_{jk}) 
       + (\delta_{ij} \delta_{kl} - \delta_{il} \delta_{jk}) 
         \left ( (\traT)^2 - \trT2 \right) \nonumber \\
       & + & 2 (\delta_{ij} T_{kl}^2 + \delta_{kl} T_{ij}^2 - \delta_{il} T_{jk}^2
       - \delta_{jk} T_{il}^2) \nonumber \\
       & + & 2 (\delta_{il} T_{jk} + \delta_{jk} T_{il} - \delta_{ij} T_{kl} 
             - \delta_{kl} T_{ij}) \traT \nonumber \\
      & = & 0.
\label{e:1.10}
\end{eqnarray}
There are (see~\cite[p.204] {Smith94}) 14 independent expressions of the form (\ref{e:1.9}). 
Applying these to
$T_{mn}$ $T_{pq}$ will yield one additional equation of the form (\ref{e:1.10}) 
which is given by:
\begin{eqnarray}
\delta^{i k m p}_{j l n q} T_{mn} T_{pq} & = &
                  2 (T_{ik} T_{jl} - T_{il} T_{jk})
                 + (\delta_{ik} \delta_{jl} - \delta_{il} \delta_{jk})
         \left ( (\traT)^2 - \trT2 \right) \nonumber \\ 
       & + & 2 (\delta_{ik} T_{jl}^2 + \delta_{jl} T_{ik}^2 - \delta_{il} T_{jk}^2
       - \delta_{jk} T_{il}^2) \nonumber \\
       & + & 2 (-\delta_{ik} T_{jl} - \delta_{jl} T_{ik} + \delta_{il} T_{jk}
             + \delta_{jk} T_{il}) \traT \nonumber \\  
      & = & 0.
\label{e:1.11}
\end{eqnarray}                            

We now require  expressions for the function $\hbfA4$ given by Eq. (\ref{e:1.8}).
The isotropic fourth--order tensor--valued function $\hbfA4$ given by Eq. (\ref{e:1.8}) is
expressible as a linear combination of 9 linear and 19 quadratic terms which are
given by:
\begin{eqnarray}
\delta_{ij} \delta_{kl} \traT, ~\delta_{ik} \delta_{jl} \traT,
~\delta_{il} \delta_{jk} \traT  \nonumber \\[3 mm]
\delta_{ij} T_{kl}, ~\delta_{ik} T_{jl}, ~\delta_{il} T_{jk}, ~\delta_{jk} T_{il},
~\delta_{jl} T_{ik}, ~\delta_{kl} T_{ij} 
\label{e:2.1}
\end{eqnarray}
and 
\begin{eqnarray}
\delta_{ij} \delta_{kl} (\traT)^2, ~\delta_{ik} \delta_{jl} (\traT)^2,
~\delta_{il} \delta_{jk} (\traT)^2, \nonumber \\[3 mm]
~\delta_{ij} \delta_{kl} \trT2, ~\delta_{ik} \delta_{jl} \trT2,
~\delta_{il} \delta_{jk} \trT2,  \nonumber \\[3 mm]
~\delta_{ij} T_{kl} \traT, ~\delta_{ik} T_{jl} \traT, ~\delta_{il} T_{jk} \traT,
~\delta_{jk} T_{il} \traT, ~\delta_{jl} T_{ik} \traT, ~\delta_{kl} T_{ij} \traT, \nonumber \\[3 mm]
~\delta_{ij} T_{kl}^2, ~\delta_{ik} T_{jl}^2, ~\delta_{il} T_{jk}^2,
~\delta_{jk} T_{il}^2, ~\delta_{jl} T_{ik}^2, ~\delta_{kl} T_{ij}^2, \nonumber \\[3 mm]
T_{ij} T_{kl} + T_{ik} T_{jl} + T_{il} T_{jk}
\label{e:2.2}
\end{eqnarray}
where
\begin{eqnarray}
\traT = T_{ii} = T_{11} + T_{22} + T_{33}, \nonumber \\[3 mm]
T_{ij}^2 = T_{ik} T_{kj}, \trT2 = T_{ij} T_{ji}.
\label{e:2.3}
\end{eqnarray}

The results given by Eq. (\ref{e:2.1}) are obtained from listing the distinct terms 
found upon applying the
15 distinct isomers of $\delta_{ij}$ $\delta_{kl}$ $\delta_{mn}$ to $T_{mn}$.
If we apply the 105 distinct isomers of 
$\delta_{ij}$ $\delta_{kl}$ $\delta_{mn}$ $\delta_{pq}$ to $T_{mn}$ $T_{pq}$,
we obtain the first 18 terms listed in (\ref{e:2.2}) together with the terms
$T_{ij}$ $T_{kl}$, $T_{ik}$ $T_{jl}$, $T_{il}$ $T_{jk}$. 
The two identities (\ref{e:1.10}) and
(\ref{e:1.11}) enable us to replace the three terms 
$T_{ij}$ $T_{kl}$, $T_{ik}$ $T_{jl}$, $T_{il}$ $T_{jk}$ by the single term
$T_{ij}$ $T_{kl}$ +  $T_{ik}$ $T_{jl}$ +  $T_{il}$ $T_{jk}$. 
It would be permissible
to employ any one of the terms $T_{ij}$ $T_{kl}$, $T_{ik}$ $T_{jl}$, $T_{il}$ $T_{jk}$
in place of $T_{ij}$ $T_{kl}$ +  $T_{ik}$ $T_{jl}$ +  $T_{il}$ $T_{jk}$ in (\ref{e:2.2}).

\section{Example of Application}

We apply the result of the previous section to obtain the expression for the fluctuating velocity, pressure--gradient correlations
as given by Eq. (\ref{e:fin}).  For convenience, and following the convention in turbulence modeling (e.g. \cite{SSG}), these
expressions are given 
in terms of the symmetric and skew--symmetric parts of
the tensor $U_{i,j}$. We employ matrix notation. Thus,
\begin{eqnarray}
U = ||U_{i,j}||,\qquad  U^T = ||U_{j,i}||,\qquad C = ||C_{ij}|| \nonumber \\
\qquad T = ||T_{ij}||,\qquad T^2 = ||T_{ik} T_{kj}||
\label{e:a.1}
\end{eqnarray}

where $U^T$ denotes the transpose of U. We have
\begin{eqnarray}
S = \frac {1}{2}(U + U^T), \qquad W = \frac{1}{2} (U - U^T), \nonumber \\
U = S + W, \qquad U^T = S - W.
\label{e:a.2}
\end{eqnarray}

 With Eqs. (\ref{e:2.1}), (\ref{e:2.2}), 
(\ref{e:a.1}) and (\ref{e:a.2}), we have:
\begin{eqnarray}
 (A_{ijkl}+A_{jikl})
\frac{  \partial U_k}{ \partial x_l} & = & a_1 S tr T + a_2 \delta_{ij} tr TS + a_3 (TS + ST) \nonumber \\
 &+ & a_4 (TW - WT)
         +  a_5 S (tr T)^2 + a_6 S tr T^2  \nonumber \\
 &+ & a_7  \delta_{ij} (tr T)^2  S +
              a_8 (TS + ST) tr T 
         +  a_9 (TW - WT) tr T  \nonumber \\
&+ & a_{10}  \delta_{ij} tr T^2S + a_{11} (T^2S + ST^2) 
         +  a_{12} (T^2W - WT^2) \nonumber \\
& + & a_{13} T (tr TS).
\label{e:a.3}
\end{eqnarray}

In obtaining Eq. (\ref{e:a.3}), it should be noted that the term $T_{ij} T_{kl} + T_{ik} T_{jl} + T_{il} T_{jk}$ appears in Eq. (\ref{e:2.2}), and
$(T_{ij} T_{kl} + T_{ik} T_{jl} + T_{il} T_{jk}) U_{k,l}$ yields
$T tr TS + 2 T S T$.  In obtaining Eq. (\ref{e:a.3}), it should also be noted that 
the following identity was employed (see~reference ~ \cite{Smith94}, p. 207):
\begin{eqnarray}
& 2 &T S T + 2 (T^2S + S T^2) - 2 (T S + S T) tr T - 2 T tr T S \nonumber \\
& - & S tr T^2 + S (tr T)^2 +2 \delta_{ij} tr T tr T S - 2 E_3 tr T^2 S = 0.
\label{e:a.4}
\end{eqnarray}

\section{Closing Remarks}
Fourth-order tensor-valued functions appear in a number of equations whose solution is required to model the behavior of a solid
undergoing deformation or a fluid in turbulent motion. The present method for representing these functions in terms of second-order tensors
that can more easily be obtained employs the  assumption that these functions are expressible as polynomials in the components of
the second-order tensors. This produces a representation that is expressible as a linear combination of a number
of linear and quadratic terms. The utility of this approach was demonstrated with respect to the fluctuating velocity, pressure-gradient
correlations of turbulence.


\end{document}